\documentclass{article}

\newtheorem{proposition}{Proposition}

\newtheorem{lemma}{Lemma}
\newtheorem{corollary}{Corollary}

\usepackage{fullpage, algorithmic, algorithm, amssymb, tikz}

\begin{document}

\title{Lexicographic identifying codes}
\author{Maximilien Gadouleau\\ School of Engineering and Computing Sciences\\ Durham University\\ m.r.gadouleau@durham.ac.uk}
\maketitle

\begin{abstract}
An identifying code in a graph is a set of vertices which intersects all the symmetric differences between pairs of neighbourhoods of vertices. Not all graphs have identifying codes; those that do are referred to as twin-free. In this paper, we design an algorithm that finds an identifying code in a twin-free graph on $n$ vertices in $O(n^3)$ binary operations, and returns a failure if the graph is not twin-free. We also determine an alternative for sparse graphs with a running time of $O(n^2d \log n)$ binary operations, where $d$ is the maximum degree. We also prove that these algorithms can return any identifying code with minimum cardinality, provided the vertices are correctly sorted.
\end{abstract}

\section{Introduction} \label{sec:intro}

Identifying codes were introduced in \cite{KCL98} for fault diagnosis in multiprocessor systems, and have since then found applications in location and detection problems. In general, an identifying code in a graph $G$ can be defined as follows. First, we denote the (closed) neighborhood of any vertex $v$ as $N(v) = \{v\} \cup \{w :  vw \in E(G)\}$. An {\em identifying code} is a subset of vertices which satisfies the following property: for any two vertices $v$ and $w$, we have $N(v) \cap C \neq N(w) \cap C \neq \emptyset$. Equivalently, it is any subset of vertices $C$ such that for all $v_1,v_2 \in V(G)$, $(N(v_1) \Delta N(v_2)) \cap C \neq \emptyset$, where $\Delta$ is the symmetric difference between two sets. A graph admits an identifying code if and only if it is {\em twin-free} \cite{CHHL07}, where twins are two vertices with the same neighborhood. We remark that the definitions above are commonly used for a so-called $1$-identifying code, where an $r$-identifying code is defined in terms of balls of radius $r$ around a vertex. Since any $r$-identifying code can be seen as a $1$-identifying code for a related graph, we do not lose any generality in considering $1$-identifying codes only. For a thorough survey of identifying codes, the reader is invited to \cite{LT08}, and an exhaustive literature bibliography on identifying codes and related topics is maintained in \cite{Lob}.

Since any superset of an identifying code is itself an identifying code, it is natural to search for the minimum cardinality $i(G)$ of an identifying code of a given graph $G$. Let us refer to an identifying code as {\em minimal} if it has no proper subset which itself is an identifying code and as {\em minimum} if it has the smallest cardinality amongst all codes. The problem of finding the minimum cardinality of an identifying code was shown to be NP-hard in \cite{LT08}. Viewing this problem as an instance of the subset cover problem \cite{CLR01}, a greedy heuristic was also designed and analyzed in \cite{LT08}. Its running time is on the order of $O(n^4)$ binary operations and has the following performance guarantees. It always finds an identifying code whose cardinality is less than $c_1 i(G) \ln n$ for some nonnegative constant $c_1$; however, there are graphs for which the algorithm always returns a code with cardinality greater than $c_2i(G) \ln n$ for another nonnegative constant $c_2$.

{\em Lexicographic codes} were introduced in \cite{Lev60} and independently rediscovered in \cite{CS86} to design large constant-weight codes, which are sets of binary vectors of equal Hamming weight with a prescribed minimum Hamming distance (see \cite{BSSS90} for a detailed review of constant-weight codes and lexicographic codes). The principle is to first sort all the vectors with the same Hamming weight, and then construct the code as we run through them. Adding a codeword is done according to a simple criterion: it must be at distance at least $d$ from the code constructed so far. The performance of the algorithm depends on the order in which the vectors have been sorted; moreover, some modifications can be added, such as starting with a predetermined set of vectors. Many record-holding constant-weight codes have been designed using lexicographic codes. However, this idea is not limited to constant-weight codes, and their application to nonrestricted binary codes has led to many interesting results \cite{Tra00}. They have also been recently applied to the construction of codes on subspaces in \cite{SE10}, also yielding record-holding codes.

In this paper, we investigate adapting the idea of lexicographic codes to identifying codes. The main contribution is an algorithm running in $O(n^3)$ binary operations which returns an identifying code for a twin-free graph, and returns a failure if the graph is not twin-free. This algorithm is then adapted to sparse graphs to run in $O(n^2d \log n)$ binary operations. Both algorithms have the same guarantees in terms of cardinality of the output. Although we are unable to give an upper bound which does not depend on the ordering of the vertices, we show that provided the vertices are properly sorted, the algorithm returns a minimum identifying code. This is fundamentally different to the greedy approach in $O(n^4)$.

\section{Algorithm for general graphs} \label{sec:main}

\subsection{Description and correctness} \label{sec:description}

Let $G$ be a graph on $n$ vertices with adjacency matrix ${\bf A}$, and let ${\bf B} = {\bf I}_n + {\bf A}$. We denote the vertices as $v_1,v_2,\ldots, v_n$, thus $b_{i,j} = 1$ if and only if $v_i \in N(v_j)$; yet we shall abuse notation and identify a vertex with its index. For instance, we refer to the vertex with minimum index in the neighborhood of $v_i$ as $\mathtt{min1}(i)$. Also, the output of our algorithm is actually the set of indices of the vertices in the code.

Before giving the pseudocode of Algorithm \ref{alg:main}, we describe it schematically below. Its input is the matrix ${\bf B}$ of the graph. It then runs along all vertices $v_j$, adding a new codeword to the code $C$ if $N(v_j) \cap C = \emptyset$ or $N(v_j) \cap C = N(v_k) \cap C$ for some $k<j$. While searching for a new codeword to add, the algorithm may return a failure if the graph is not twin-free, which we identify as $n+1 \in C$. After the $j$-th step, the code $C$ then `identifies' the first $j$ vertices, i.e. they are all covered in a distinct fashion. We keep track of the intersections $N(v_i) \cap C$ in a matrix ${\bf X}$. After going through all vertices, the algorithm then returns an identifying code $C$ or a failure (if $n+1 \in C$) if the graph is not twin-free.

\begin{algorithm}

\caption{Main algorithm for general graphs}
\label{alg:main}

\begin{algorithmic}
\STATE $C \gets \emptyset$, $X \gets {\bf 0}_n$, $j \gets 1$
\WHILE{$j\leq n$ and $n+1 \notin C$}
	\STATE $l \gets 0$
    \IF[\ensuremath{v_j} is not covered]{${\bf X}(j) = {\bf 0}$}
        \STATE $l \gets \mathtt{min1}(j)$
    \ELSE
        \STATE $k \gets 1$
        \WHILE[\ensuremath{v_j} is covered, so we search if it is identified]{${\bf X}(j) \neq {\bf X}(k)$ and $k<j$}
            \STATE $k \gets k+1$
        \ENDWHILE
%	 	  \STATE $k \gets \mathtt{search}(j,{\bf X})$
        \IF[\ensuremath{v_j} is not identified]{$k < j$}
            \STATE $l \gets \mathtt{min2}(j,k)$ %\COMMENT{add a new codeword to identify $v_j$}
        \ENDIF
    \ENDIF
    \IF[A new codeword has been found]{$1 \leq l \leq n$}
        \STATE $C \gets C \cup \{l\}$ %\COMMENT{update the code}
        \STATE ${\bf X}^T(l) \gets {\bf B}^T(l)$ %\COMMENT{update the code matrix}
    \ENDIF
    \STATE $j \gets j+1$
\ENDWHILE
\RETURN $C$
\end{algorithmic}
\end{algorithm}

The subroutine $\mathtt{min2}(j,k)$ returns the first vertex which identifies $v_j$ if it exists and a failure otherwise, i.e. it determines the first vertex in lexicographic order in $N(v_j) \Delta N(v_k)$. If $N(v_j) = N(v_k)$, then it returns $n+1$. It is given in Algorithm \ref{alg:newcodeword_main}.

%Input: two integers $j$ and $k$ between $1$ and $n$\\
%Output: $l$, where $v_l$ is the first vertex in lexicographic order such that $\{v_k,v_l\}$ identify $v_j$ if such exists, a failure ($l=n+1$) otherwise\\

\begin{algorithm}

\caption{$\mathtt{min2}(j,k)$ subroutine}
\label{alg:newcodeword_main}

\begin{algorithmic}
\STATE $l \gets 1$
\WHILE{$l \leq n$ and ${\bf B}(j,l) = {\bf B}(k,l)$}
    \STATE $l \gets l+1$
\ENDWHILE
%\IF{$l = n+1$}
%    \STATE $f \gets 1$
%\ENDIF
\RETURN $l$
\end{algorithmic}
\end{algorithm}

%\begin{algorithm}
%\caption{$\mathtt{search}(j,K)$}
%\label{alg:search}
%\begin{algorithmic}
%\STATE $i \gets K$, $k \gets 1$
%\WHILE{$i\geq 1$ and $k < j$}
%	\WHILE{$i \geq 1$ and ${\bf X}(k,i) = {\bf X}(j,i)$}
%		\STATE $i = i-1$
%	\ENDWHILE
%	\STATE $k = k+1$
%\ENDWHILE
%\RETURN $k$
%\end{algorithmic}
%\end{algorithm}

We now justify this claim in Lemma \ref{lemma:subroutine} below.

\begin{lemma} \label{lemma:subroutine}
The subroutine $\mathtt{min2}(j,k)$ returns the minimum element in $N(v_j) \Delta N(v_k)$ if this symmetric difference is non-empty, and a failure ($l=n+1$) otherwise.
\end{lemma}

\textbf{Proof}
First, if $N(v_j) = N(v_k)$, then ${\bf B}(j,l) = {\bf B}(k,l)$ for all $1 \leq l \leq n$. Therefore, the \textbf{while} loop will only stop once $l=n+1$, and hence the subroutine returns a failure. Second, if $N(v_j) \neq N(v_k)$, then the minimum element in $N(v_j) \Delta N(v_k)$ is the smallest $l$ such that ${\bf B}(j,l) \neq {\bf B}(k,l)$. It is clear that the subroutine returns this value.
$\Box$

\begin{proposition} \label{prop:correctness}
Algorithm \ref{alg:main} returns an identifying code if the input graph is twin-free, and a failure ($n+1 \in C$) otherwise.
\end{proposition}

\textbf{Proof}
First of all, we prove that the algorithm returns a failure if and only if the graph is not twin-free. In the latter case, let $k$ be the smallest integer such that the set $\{i \neq k : N(v_k) = N(v_i) \}$ is not empty, and let $j$ be the minimum element of this set (hence $k < j$, $N(v_k) = N(v_j)$). It is easily shown that after the $k$-th step, $v_k$ is covered. On the $j$-th step, Algorithm \ref{alg:main} first checks if $v_j$ is covered. Since $v_k$ is covered and $N(v_k) = N(v_j)$, then $v_j$ is also covered. Algorithm \ref{alg:main} then finds that $k$ is the smallest integer satisfying ${\bf X}(k) = {\bf X}(j)$, and hence calls the subroutine $\mathtt{min2}(j,k)$. By Lemma \ref{lemma:subroutine} this returns a failure, and hence the whole algorithm returns a failure. Conversely, the only case where the subroutine (and hence the algorithm) returns a failure is when there exist $k < j$ such that $N(v_j) = N(v_k)$, i.e. the graph is not twin-free.

We now assume that the graph is twin-free, and hence we have $l \leq n$ at any step. We need to show that the output $C$ of Algorithm \ref{alg:main} is an identifying code. Let us denote the matrix ${\bf X}$ and the code $C$ obtained after $j$ steps as ${\bf X}^j$ as $C^j$, respectively. Note that for all $a$, ${\bf X}^j(a)$ reflects how the vertex $v_a$ is covered by $C^j$: $N(v_a) \cap C^j = \mathrm{supp}({\bf X}(a)) = \{b : {\bf X}^j(a,b) = 1\}$. The following claim is the cornerstone of the proof.

\textbf{Claim}: After step $j$, all ${\bf X}^j(i)$'s are nonzero and distinct for $1 \leq i \leq j$. 

The proof goes by induction on $j$, and is trivial for $j=1$. Suppose it is true for $j-1$, then 
\begin{equation} \nonumber
	\mathrm{supp}({\bf X}^{j-1}(a)) = N(v_a) \cap C^{j-1} \subseteq N(v_a) \cap C^j = \mathrm{supp}({\bf X}^j(a)).
\end{equation}
It is hence easy to show that if ${\bf X}^{j-1}(a) \neq {\bf 0}$, then ${\bf X}^j(a) \neq {\bf 0}$ and if ${\bf X}^{j-1}(a) \neq {\bf X}^{j-1}(b)$, then ${\bf X}^j(a) \neq {\bf X}^j(b)$ for all $a$ and $b$. It immediately follows that the vectors ${\bf X}^j(i)$'s are all nonzero and distinct for $1 \leq i \leq j-1$, and we only have to consider ${\bf X}^j(j)$. Three cases occur when the algorithm reaches step $j$.
\begin{itemize}
    \item Case I: ${\bf X}^{j-1}(j)$ is nonzero and distinct to any ${\bf X}^{j-1}(i)$ for $1 \leq i \leq j-1$. Then as shown above, ${\bf X}^j(j)$ is nonzero and distinct to all ${\bf X}^j(i)$'s.

    \item Case II: ${\bf X}^{j-1}(j)$ is nonzero and equal to ${\bf X}^{j-1}(k)$ for some $k < j$. First, we remark that $k$ is unique, as ${\bf X}^{j-1}(k) \neq {\bf X}^{j-1}(i)$ for all other $i$. The $\mathtt{min2}(k,j)$ subroutine then returns an element $v_l \in N(v_j) \Delta N(v_k)$, and hence ${\bf X}^j(j,l) \neq {\bf X}^j(k,l)$.

    \item Case III: ${\bf X}^{j-1}(j) = {\bf 0}$. Then by hypothesis ${\bf X}^{j-1}(j) \neq {\bf X}^{j-1}(i)$ for all $1 \leq i \leq j-1$, and hence ${\bf X}^j(j) \neq {\bf X}^j(i)$. Also, ${\bf X}^j(j)$ is the unit vector ${\bf e}_{\mathtt{min1}(j)}$, which is nonzero.
\end{itemize}

Therefore, for the code $C^n = C$ obtained after $n$ steps, $N(v_a) \cap C$ are all nonzero and distinct for all $1 \leq a \leq n$. It is hence an identifying code.
$\Box$

\subsection{Performance} \label{sec:performance}

We now investigate the performance of Algorithm \ref{alg:main}. We are first interested in the cardinality of its output. Clearly, this depends on the order in which the vertices are sorted. We show below that provided the order is suitable, the algorithm can find any minimal identifying code, and hence can return a minimum one.

\begin{proposition} \label{prop:minimal}
Suppose that the graph is twin-free and that $M = \{v_1,v_2,\ldots,v_m\}$ forms an identifying code. Then Algorithm \ref{alg:main} returns an identifying code that is a subset of $M$.
\end{proposition}

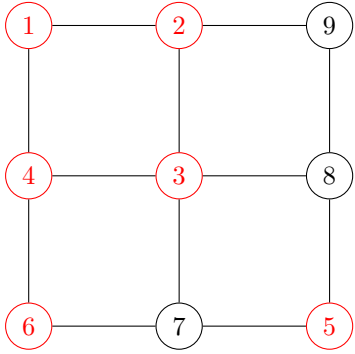
\begin{figure}
\centering
\begin{tikzpicture}
	\tikzstyle{every node}=[draw,shape=circle];
	
	\node[color=red] (1) at (0,4) {$1$};
	\node[color=red] (2) at (2,4) {$2$};
	\node (3) at (4,4) {$9$};
	\node[color=red] (4) at (0,2) {$4$};
	\node[color=red] (5) at (2,2) {$3$};
	\node (6) at (4,2) {$8$};
	\node[color=red] (7) at (0,0) {$6$};
	\node (8) at (2,0) {$7$};
	\node[color=red] (9) at (4,0) {$5$};

	\draw (1)--(2);
	\draw (2)--(3);
	\draw (4)--(5);
	\draw (5)--(6);
	\draw (7)--(8);
	\draw (8)--(9);
	
	\draw (1)--(4);
	\draw (4)--(7);
	\draw (2)--(5);
	\draw (5)--(8);
	\draw (3)--(6);
	\draw (6)--(9);

\end{tikzpicture}
\caption{A graph and a sorting of vertices such that the lexicographic code is not minimal}
\label{fig:example}
\end{figure}

\textbf{Proof}
We know by Proposition \ref{prop:correctness} that the algorithm returns an identifying code; we only have to prove that all codewords are in $M$. At step $j$, three cases need to be distinguished.
\begin{itemize}
    \item Case I: $v_j$ is covered and identified, then no codeword is added.

    \item Case II: $v_j$ is covered but not identified, i.e. $(N(v_j) \Delta N(v_k)) \cap C^{j-1} = \emptyset$ for some $k < j$. The subroutine returns the smallest element $v_l$ in $N(v_j) \Delta N(v_k)$. Since $M$ is an identifying code, the set $(N(v_j) \Delta N(v_k)) \cap M$ is not empty, hence $v_l \in M$.

    \item Case III: $v_j$ is not covered. The algorithm then selects the next codeword to be $\mathtt{min1}(j)$, which is necessarily in $M$ as $N(v_j) \cap M \neq \emptyset$.
\end{itemize}
Therefore, the algorithm only adds codewords of $M$, and hence returns a subcode of $M$.
$\Box$
We remark that Algorithm \ref{alg:main} does not necessarily return a minimal code, as seen in Figure \ref{fig:example}. Algorithm \ref{alg:main} would return the code $\{1,2,3,4,5,6\}$ while $\{2,3,4,5,6\}$ is a minimal identifying code.

On the other hand, if $M$ is minimal, then it has no proper subset that itself is an identifying code; Algorithm \ref{alg:main} thus returns it. We obtain the following corollary.

\begin{corollary} \label{cor:minimal}
Provided that the vertices are sorted such that $v_1,v_2,\ldots,v_m$ form a minimal identifying code for some $1 \leq m \leq n$, Algorithm \ref{alg:main} will return this identifying code.
\end{corollary}

%Therefore, our algorithm can return a minimum identifying code with probability at least $\frac{1}{{n \choose K}}$, where $K$ is the cardinality. Furthermore, for all $0 \leq x \leq n-K$, Algorithm \ref{alg:main} returns an identifying code with cardinality at most $K+x$ with probability at least $\frac{{n-K \choose x}}{{n \choose K+x}}$. GIVE PROOF

Proposition \ref{prop:minimal} also implies that the probability that the output has cardinality no more than $K$ is at least the probability that the first $K$ vertices form an identifying code. Hence
 our algorithm returns a minimum identifying code with probability at least $\frac{1}{{n \choose i(G)}}$.

\begin{proposition} \label{prop:running_time_main}
The running time of Algorithm \ref{alg:main} is $O(n^3)$ binary operations.
\end{proposition}

\textbf{Proof}
Clearly, we have to run the iteration for $j$ exactly $n$ times. For each iteration, the step demanding the highest number of operations is the search for $k$. We consider at most $j-1$ values of $k$, comparing at most $n$ bits to verify whether ${\bf X}(j) \neq {\bf X}(k)$. Therefore, the running time is $O(n^3)$.
$\Box$

\section{Algorithm for sparse graphs} \label{sec:sparse}

For sparse graphs, it is more efficient not to work with the whole adjacency matrix, but with the neighborhood array $A \in \mathcal{P}(E)^n$, defined as $A(v_i) = N(v_i)$, where the neighborhood is sorted in increasing lexicographic order. Then, instead of adding the column of the adjacency matrix corresponding to a new codeword, we only update the code array $X(v)$ for all vertices adjacent to the new codeword. THe algorithm for sparse graphs is given in Algorithm \ref{alg:sparse}; its input is the neighorhood array, and it returns an identifying code $C$ or a failure ($n+1 \in C$) if the graph is not twin-free.

\begin{algorithm}

\caption{Main algorithm for sparse graphs}
\label{alg:sparse}

\begin{algorithmic}
\STATE $C \gets \emptyset$, $X \gets \emptyset^n$, $j \gets 1$, $f \gets 0$
\WHILE{$j\leq n$ and $n+1 \notin C$}
	\STATE $l \gets 0$
    \IF[$v_j$ not covered]{$X(j) = \emptyset$}
        \STATE $l \gets A(j,1)$
    \ELSE
        \STATE $m \gets X(j,1)$, $k \gets 1$
        \WHILE{$X(j) \neq X(k)$ and $k<j$}
            \STATE $k \gets k+1$
        \ENDWHILE
        \IF[$v_j$ not identified]{$k < j$}
            \STATE $l \gets \mathtt{min3}(j,k)$
        \ENDIF
    \ENDIF
    \IF{$1 \leq l \leq n$}
        \STATE $C \gets C \cup \{l\}$
        \FOR{$i$ from $1$ to $d_l$}
            \STATE $X(A(l,i)) \gets X(A(l,i)) \cup \{l\}$
        \ENDFOR
    \ENDIF
    \STATE $j \gets j+1$
\ENDWHILE
\RETURN $C$
\end{algorithmic}
\end{algorithm}

Similar to the general case, the $\mathtt{min3}(j,k)$ subroutine produces the first vertex $v_l$ which identifies $v_j$ if it exists and a failure otherwise, i.e. it determines the first vertex in lexicographic order which covers either $j$ or $k$, but not both. It is given in Algorithm \ref{alg:newcodeword}.
%Input: two vertices $j$ and $k$\\
%Output: $l$, where $v_l$ is the first vertex in lexicographic order such that $\{v_k,v_l\}$ identify $v_j$ if such exists, a failure otherwise\\

\begin{algorithm}

\caption{$\mathtt{min3}(j,k)$ subroutine}
\label{alg:newcodeword}

\begin{algorithmic}
\STATE $l \gets n+1$
\WHILE{$a \leq \min\{d_j,d_k\}$}
    \IF{$A(j,a) \neq A(k,a)$}
        \STATE $l \gets \min\{A(j,a),A(k,a)\}$
    \ENDIF
    \STATE $a \gets a+1$
\ENDWHILE
\IF{$l=n+1$}
    \IF{$d_j < d_k$}
        \STATE $l \gets A(k,d_j+1)$
    \ELSIF{$d_k < d_j$}
        \STATE $l \gets A(j,d_k+1)$
    \ENDIF
\ENDIF

\RETURN $l$
\end{algorithmic}
\end{algorithm}

The same results on correctness and the possibility of returning a minimum code also hold for Algorithm \ref{alg:sparse}; they are summarized below.

\begin{proposition} \label{prop:sparse}
If the graph is not twin-free, then Algorithm \ref{alg:sparse} returns a failure. Otherwise, the algorithm returns an identifying code contained in $\{v_1,v_2,\ldots,v_m\}$, where $m$ is the minimum integer such that this forms an identifying code.

The running time of Algorithm \ref{alg:sparse} is $O(n^2d \log n)$ binary operations.
\end{proposition}

\textbf{Proof}
The proof of correctness of Algorithm \ref{alg:sparse} is similar to that of Algorithm \ref{alg:main}, and is hence omitted. We hence determine the running time of the algorithm.
$\Box$

\bibliographystyle{IEEEtran}
\bibliography{g}

\end{document}